\documentclass[preprint,11pt,3p,review,sort&compress]{elsarticle}
\usepackage{amsmath}
\usepackage{amsfonts}
\usepackage{amssymb}
\usepackage{xcolor}
\definecolor{darkGreen1}{rgb}{0,0.5,0}
\definecolor{violet1}{rgb}{0.26,0,0.53}
\definecolor{orange1}{rgb}{0.85,0.33,0.1}

\usepackage{tikz}
\usetikzlibrary{calc,shapes,arrows,patterns,decorations.pathmorphing,decorations.markings}
\usetikzlibrary{arrows.meta}
\usetikzlibrary{external}
\usepackage{bm}

\begin{document}
	\title{$H_\infty$ optimization of multiple tuned mass dampers for multimodal vibration control}
	\author[s3l]{G.~Raze\corref{cor1}}
	\ead{g.raze@uliege.be}
	\author[s3l]{G.~Kerschen}
	\ead{g.kerschen@uliege.be}
	
	\address[s3l]{Space Structures and Systems Laboratory, Aerospace and Mechanical Engineering Department, University of Li\`{e}ge \\
Quartier Polytech 1 (B52/3), All\'{e}e de la D\'{e}couverte 9, B-4000 Li\`{e}ge, Belgium}
	\cortext[cor1]{Corresponding author}
	

	\begin{abstract}
	    In this paper, a new computational method for the purpose of multimodal vibration mitigation using multiple tuned mass dampers is proposed. Classically, the minimization of the maximum amplitude is carried out using direct $H_\infty$ optimization. However, as shall be shown in the paper, this approach is prone to being trapped in local minima, in view of the nonsmooth character of the problem at hand. This is why this paper presents an original alternative to this approach through norm-homotopy optimization. This approach, combined with an efficient technique to compute the structural response, is shown to outperform direct $H_\infty$ optimization in terms of speed and performance. Essentially, the outcome of the algorithm leads to the concept of \textit{all-equal-peak design} for which all the controlled peaks are equal in amplitude. This unique design is new with respect to the existing body of knowledge.
	\end{abstract}
	

	\maketitle

    \section{Introduction}
    
		Tall, slender and light structures are more and more used in various engineering fields for performance, compliance with regulations and/or esthetic reasons. An inconvenient feature of these structures is their proneness to exhibit lightly-damped, high-amplitude resonances. Such resonances may shorten the lifetime of these structures and even render them dangerous for human use. A possible cure against this is to attach a tuned mass damper (TMD) to the structure in the attempt to mitigate its vibratory amplitude.
		
		The TMD was originally proposed by Frahm~\cite{frahm1911device} as a spring-mass device for suppressing a specific resonance frequency. Ormondroyd and Den Hartog~\cite{DenHartog1928} added a damper to this device and tuned the resulting absorber based on fixed points of the compliance (i.e., the transfer function between a displacement of interest and the external forcing amplitude) of the controlled structure. The fixed points, independent on absorber damping, were chosen to be equal in amplitude. Brock~\cite{brock1946note} proposed a particular value of the absorber damping coefficient such that the two fixed points simultaneously be approximate maxima of the compliance. Those developments, gathered in Den Hartog's book~\cite{den1985mechanical}, laid down the foundations of the {\it equal-peak design} because the controlled compliance exhibits two peaks of (approximately) equal amplitude, usually much lower than that of the uncontrolled structure. Since then, there has been a large number of  tuning formulas varying with the loading conditions and objectives at hand. For instance, Warburton~\cite{Warburton1982} who proposed an unified approach for the fixed-point method gave optimal parameters for several cases. Nishihara and Asami~\cite{nishihara2002closed} found the exact analytical solution to the $H_\infty$-optimization problem by minimizing the maximum value of the compliance under the assumption that the latter exhibits two peaks of equal amplitude (i.e., an exact equal-peak design). 
		
		TMDs are used in a wide range of civil and mechanical engineering applications. A review on the subject has been done recently by Elias and Matsagar~\cite{Elias2017}. However, being \textit{tuned} to a particular frequency, the TMD may feature a lack of robustness when the targeted resonance frequency is uncertain or varies with time. A recent approach proposed by Dell'Elce et al~\cite{DellElce2018} tunes the absorber parameters according to the maximum uncertainty on the host structure. Structural nonlinearities may also detune the absorber, but their effect can be countered effectively using a nonlinear tuned vibration absorber~\cite{Habib2015a}. Alternatively, a number of small TMDs tuned over a frequency band centered around the resonance frequency of interest can be robust to variations in that frequency, but also more efficient than a single TMD in the sense that it yields a smaller minimum of maximum amplification. The beneficial effects of a TMD array were first discovered by Snowdon ~\cite{Snowdon1964} and Iwanami and Seto~\cite{Iwanami1984}, but the true potential of multiple tuned mass dampers (MTMDs) was unlocked in the works of Igusa and Xu~\cite{igusa1994vibration}, Yamaguchi and Harpornchai~\cite{Yamaguchi1993} and Abe and Fujino~\cite{Abe1994}, among others.
		
		MTMD can also target multiple resonances by assigning one or several TMDs per mode to be controlled. Early works about multimodal vibration mitigation used bars~\cite{Neubert1964} and beams~ \cite{Snowdon1966,Kitis1983,Ozguven1986} as host structures. Rana and Soong~\cite{Rana1998} applied this approach to spring-mass systems and, as their discussion reveals, this second use of MTMD received less attention than the first one. In those studies, the absorbers were tuned such that the controlled compliance displays two pairs of equal peaks in place of the first two resonances. Clark~\cite{Clark1988} demonstrated the MTMD efficiency in reducing the maximum acceleration experienced at the top of a building during an earthquake. Yau and Yang~\cite{Yau2004} robustly controlled two modes of cable-stayed bridges traveled by high-speed trains, by using one TMD array per mode to be controlled.
		
		Closed-form expressions for the absorber parameters are usually available when the absorber is placed on undamped single-degree-of-freedom oscillators. They can also be used for multiple-degree-of-freedom structures provided that their resonance frequencies are widely spaced. Real-life structures always violate these assumptions to some extent. Krenk and H{\o}gsberg~\cite{Krenk2016} proposed to use quasi-static and quasi-dynamic background correction terms to account for non-resonant modes. Several numerical optimization techniques were used to tune TMD and MTMD parameters; examples include parameters space exploration~\cite{Li2006,Han2008}, gradient-based optimization \cite{Warburton1982,Rade2000,Zuo2004,Hoang2005,Li2007,Piccirillo2019}, metaheuristic optimization (such as particle swarm optimization~\cite{Leung2008,Leung2009}, genetic algorithms~\cite{Hadi1998,Arfiadi2011,Mohebbi2012}, harmony search~\cite{Bekdas2011,Nigdeli2017,Zhang2017}, ant colony optimization~\cite{Viana2008,Farshidianfar2013}, simulated annealing~\cite{Yuan2017} and coral reefs optimization~\cite{Salcedo-Sanz2017}) and hybrid optimization algorithms, using both metaheuristic and gradient-based optimizations~\cite{Yang2014}. The aforementioned optimization procedures suffer from at least one of the following limitations. First, some of them neglect the effect of damping and/or other non-resonant structural modes in the structure. Second, metaheuristic optimization algorithms can be prohibitive in terms of computational cost when the number of variables to optimize becomes large. Finally, one or several absorber parameters are often assumed to be fixed, which may lead to a suboptimal design.
		
		This paper proposes a novel MTMD tuning methodology for multimodal vibration mitigation of linear structures. The developed algorithm finds the resonance peaks of the compliance and minimizes their amplitude simultaneously. It results in a so-called {\it all-equal-peak design}, i.e., all the peaks of the controlled resonances are equal in amplitude. The paper is organized as follows. In Section~\ref{sec:allequalpeaks}, the general principles of the tuning methodology approach and of the optimization algorithm are introduced. Section~\ref{sec:optimization} details the numerical optimization procedure. Section~\ref{sec:examples} then illustrates the concept of all-equal-peak design with a simple spring-mass system and a simply supported plate featuring high modal density. Finally, the conclusions of the present study are drawn in Section~\ref{sec:conclusion}.
        
    \section{A norm-homotopy approach for $H_\infty$ optimization}
    \label{sec:allequalpeaks}
    
    	In this paper, the structure without absorbers and with $N_a$ attached absorbers is referred to as \textit{host structure} and \textit{controlled structure}, respectively. In the presence of harmonic forcing, the vibratory amplitude of a single-degree-of-freedom host structure is classically mitigated through the minimization of the $H_\infty$ norm of a given transfer function, i.e., its maximum amplitude, resulting in the equal-peak design for which there exist well-established analytical tuning rules, e.g., \cite{nishihara2002closed}. However, multiple-degree-of-freedom host structures have more complicated transfer functions, which rules out the possibility of tuning the absorbers analytically. Resorting to numerical optimization for minimizing the $H_\infty$ norm is then necessary.
    	
\begin{figure}[!ht]
			\centering
			\includegraphics[width=0.7\textwidth]{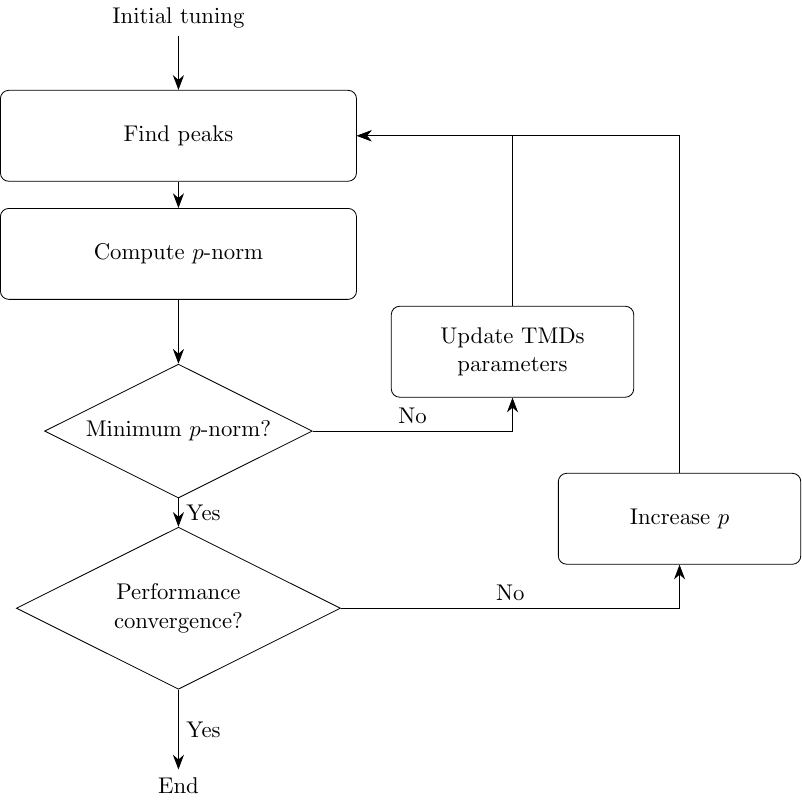}
			\caption{Conceptual flowchart of the proposed norm-homotopy optimization algorithm.}
			\label{fig:flowchart0}
		\end{figure}

    	When considering multiple resonances, one inherent difficulty with the $H_\infty$ norm is that it considers only the resonance peak exhibiting the largest amplitude, i.e., it disregards the other controlled peaks. Their amplitude is thus minimized later in the optimization process when they themselves feature the largest amplitude. This typically results in a nonsmooth cost function which may lead to a premature termination of the algorithm. The alternative strategy proposed in this paper relies on a norm-homotopy optimization during which problems of increasing complexity are solved sequentially using the previously-obtained parameters as an initial guess for the next problem. Specifically, the $p$-norm of the vector containing the controlled peak amplitudes, i.e., $\|x\|_p=\left(\sum_{i=1}^n|x_i|^p\right)^{1/p}$, is minimized, and $p$ is sequentially increased so as to approach the $H_\infty$ norm, as schematically presented in Fig.~\ref{fig:flowchart0}. A low value of $p$ puts more weight on the resonance peaks with lower amplitudes and makes the optimization problem less stiff, whereas the subsequent increase in $p$ ensures that resonances with large amplitudes are penalized enough.
    	
\begin{figure}[!ht]
		    \centering
				\includegraphics[width=0.6\textwidth]{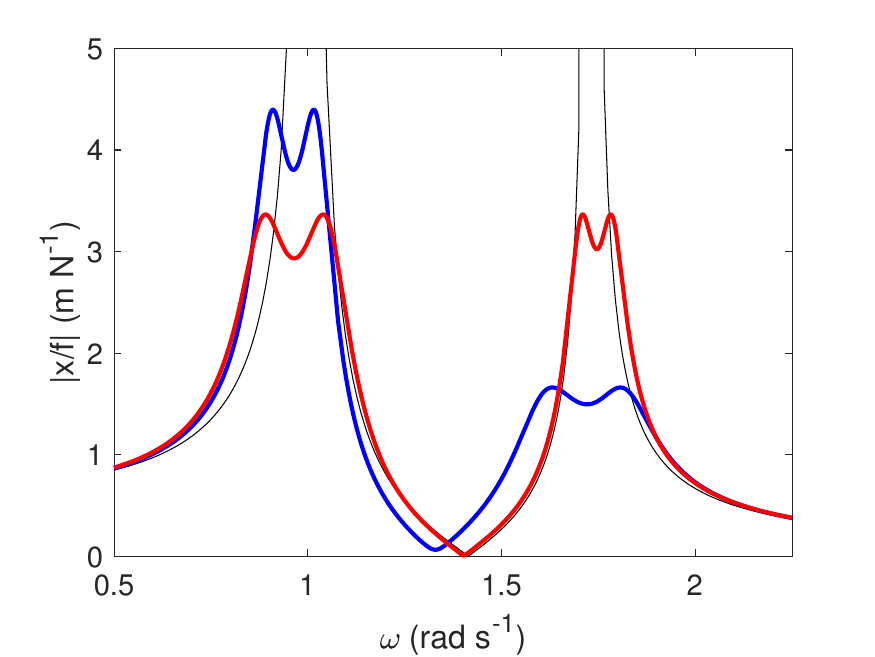}
				\caption{MTMD applied to a two-degree-of-freedom system: uncontrolled structure (\textbf{---}), solution from \cite{Ozguven1986} (\textcolor{blue}{\textbf{---}}) and norm-homotopy solution (\textcolor{red}{\textbf{---}}).}
				\label{fig:OC_Comparison}
		\end{figure}

		A typical output of the norm-homotopy optimization is shown in Fig.~\ref{fig:OC_Comparison} for the mitigation of the resonances of a two-degree-of-freedom system (studied more in depth in Section~\ref{ssec:2dofs}). Clearly, the algorithm is able to enforce the same amplitude for the four resonances. Existing algorithms in the literature, see, e.g., \cite{Ozguven1986}, can also enforce equal peaks for each resonance, but the amplitudes associated with each pair of peaks are not equa, and the transfer function thus exhibits a higher $H_\infty$ norm. This all-equal-peak design appears as a generalization to multiple modes of the equal-peak design, and can only be achieved through numerical optimization, given the complexity of the problem at hand.

	\section{The proposed optimization algorithm}
	\label{sec:optimization}
		
\begin{figure}[!ht]
					\centering
					\includegraphics[width=0.7\textwidth]{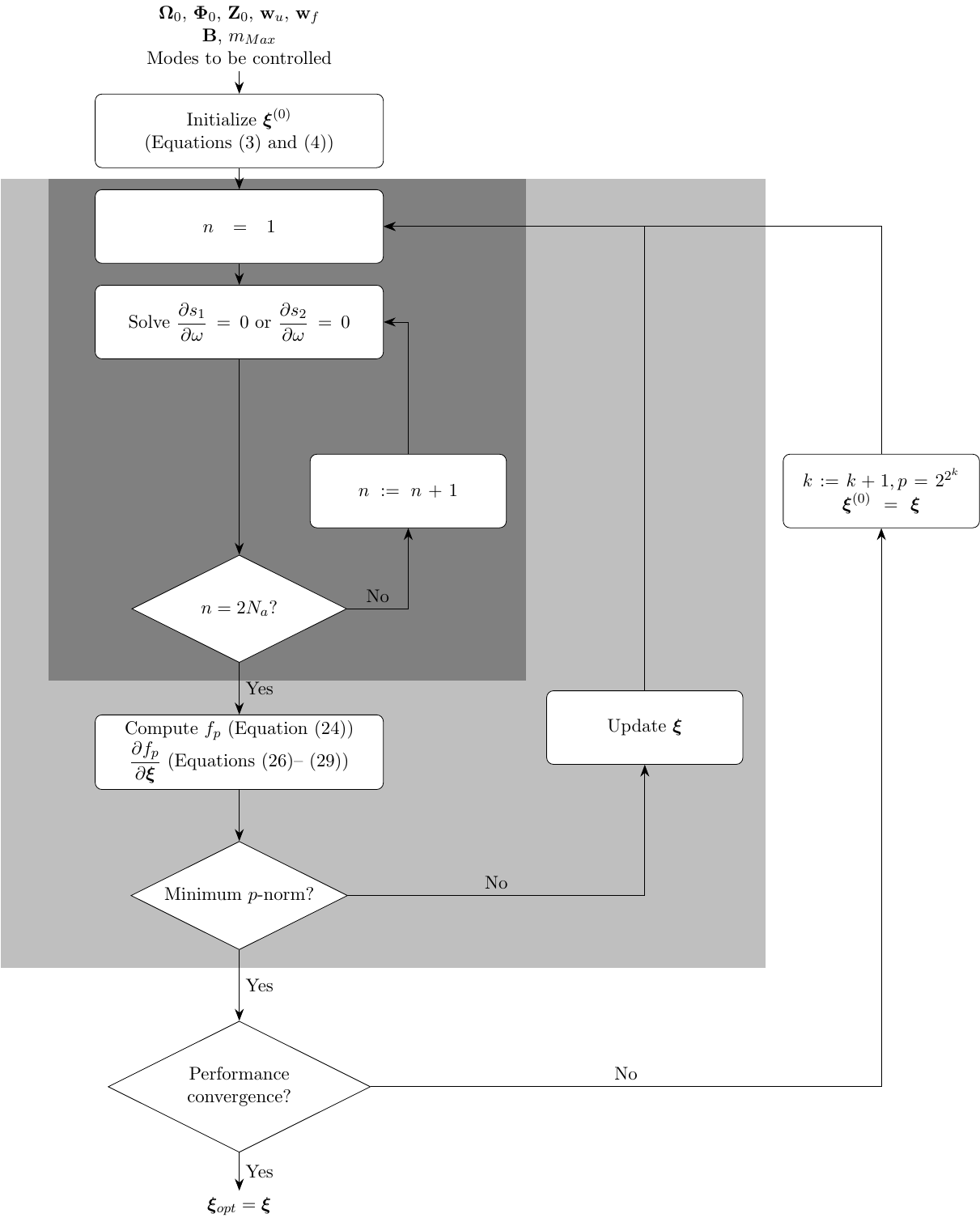}
					\caption{Flowchart of the norm-homotopy optimization algorithm; the dark gray area represents the peak-finding algorithm, and the pale gray area represents the $p$-norm optimization algorithm.}
					\label{fig:flowchart}
				\end{figure}

		The proposed algorithm is presented in Fig.~\ref{fig:flowchart}. The objective is to optimize the parameters of the different absorbers gathered in a vector $\bm{\xi}$. A tuning based on the well-established single-degree-of-freedom formulas from the literature provides an initial guess $\bm{\xi}_0$, as discussed in Section \ref{ssec:initialTuning}. Because damping in the host structure and non-resonant modes are ignored and because there is cross-influence between the absorbers, the resulting performance is usually not satisfactory, and the parameters have to be further optimized. Once the transfer function is computed for a specific set of parameters, the resonance peaks have to be located. The strategy for peak finding is described in Subection~\ref{ssec:peaks}, whereas the optimization procedure of the different cost functions is presented in Subsections~\ref{ssec:pnorm} and ~\ref{norm}. We note that the algorithm relies extensively on evaluations of the compliance, which may become computationally expensive for structures with a large number of degrees of freedom. To cope with this issue, Section~\ref{ssec:dynamics} formulates the dynamics of the controlled structure using the Sherman-Morrison-Woodbury formula~\cite{woodbury1950inverting}.
		 
			\subsection{Initial tuning}		
			\label{ssec:initialTuning}
			
			A MTMD with one TMD per targeted mode is considered herein, and we assume that the $n^{th}$ absorber targets the $r^{th}$ resonant mode of the host structure. It may be shown (see~\ref{anx:TMD} or~\cite{Warburton1982}) that, if one neglects non-resonant modes, the structure acts from the absorber point of view as an equivalent one-degree-of-freedom mechanical oscillator with the following modal mass, damping and stiffness: 
    		\begin{equation}
    			m_r = \frac{1}{\phi_{a,n}^2}, \qquad c_r = 2\omega_r m_r \zeta_r, \qquad k_r = \omega_r^2 m_r
    			\label{eq:modal1dofCharac}
    		\end{equation}
    		where $\phi_{a,n}$ is the $r^{th}$ mass-normalized mode shape of the host structure at the location where the absorber is attached, $\omega_r$ is the resonance frequency, and $\zeta_r$ is the modal damping ratio. Classical formulas from the literature can then be used for absorber tuning (\cite{den1985mechanical,Warburton1982,nishihara2002closed} or even \cite{Asami2002} if damping in the host structure is taken into account). In this paper, the formulas from Nishihara and Asami~\cite{nishihara2002closed} are used. Since they are exact in the single-degree-of-freedom case, they are expected to be a reasonably accurate initial guess for a multiple-degree-of-freedom case. From Eq.~\eqref{eq:modal1dofCharac}, the modal mass ratio is defined as
    		\begin{equation}
    			\mu_{a,n} = \frac{m_{a,n}}{m_r},
    			\label{eq:mmr}
    		\end{equation}
    		and the absorber stiffness and damping are computed as
    		\begin{equation}
    			k_{a,n} = \frac{8}{\left(1+\mu_{a,n}\right)^2}\frac{16+23\mu_{a,n}+9\mu_{a,n}^2+2(2+\mu_{a,n})\sqrt{4+3\mu_{a,n}}}{3(64+80\mu_{a,n}+27\mu_{a,n}^2)} \omega_r^2 m_{a,n}
    			\label{eq:stiffnessInit}
    		\end{equation}
    		\begin{equation}
    			c_{a,n} = \frac{1}{2}\sqrt{\frac{8+9\mu_{a,n}-4\sqrt{4+3\mu_{a,n}}}{1+\mu_{a,n}}}\sqrt{k_{a,n}m_{a,n}} ,
				\label{eq:dampingInit}
    		\end{equation}
    		respectively. With these formulas, the maxima of the compliance are expected to be near
			\begin{equation}
				\omega_{n1,n2} = \dfrac{1}{1+\mu_{a,n}}\left(1\pm\sqrt{\dfrac{\mu_{a,n}}{2+\mu_{a,n}}}\right)\omega_r.
				\label{eq:peaksInit}
			\end{equation}			    		
    		If maximum efficiency is sought, the modal mass ratio should be maximized~\cite{nishihara2002closed}, which, according to Eq.~\eqref{eq:mmr}, is equivalent to minimizing the modal mass $m_r$. Going back to Eq.~\eqref{eq:modal1dofCharac}, the modal mass is minimized if the absorber is placed at a maximum of modal amplitude of the $r^{th}$ mode in the host structure. This result is by no means new; further considerations are given in Petit et al~\cite{Petit2009} when either this location is not acceptable for attaching an absorber, or when the activity of the neighboring modes is too prominent.
    		    		
\begin{figure}[!ht]
				\centering
				\includegraphics[width=0.6\textwidth]{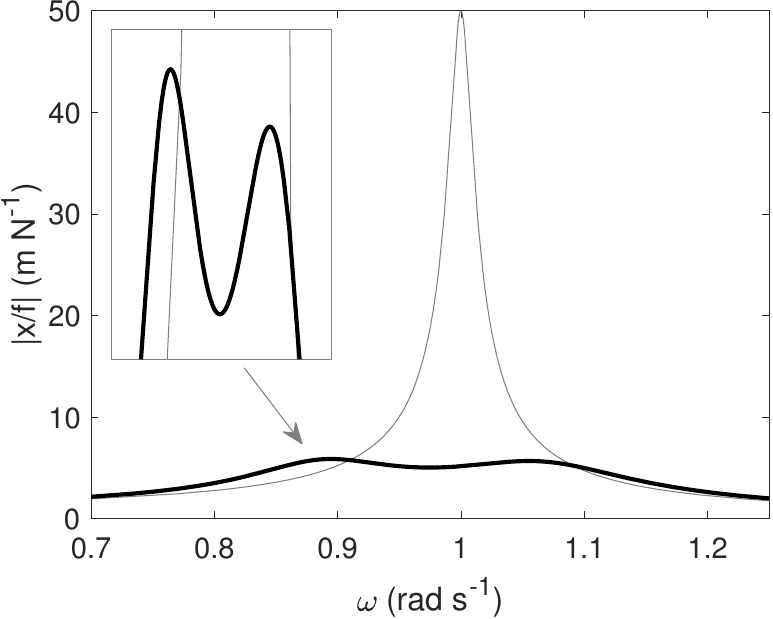}
				\caption{Compliance of the damped single-degree-of-freedom system: uncontrolled (\textcolor{gray}{---}) and controlled (\textbf{---}) structure.}
				\label{fig:sdofDampedInit}
			\end{figure}

    		This procedure can be repeated for each absorber to yield an initial design for the attached MTMD. Because damping in the host structure and non-resonant modes are ignored and because there is cross-influence between the absorbers, the absorber parameters usually have to be further optimized. For illustration, a damped single-degree-of-freedom host structure controlled by a single TMD, is studied through Section~\ref{sec:optimization}. Its purpose is to demonstrate the working principles of the algorithm with a simple example. The parameters of the host structure are $m_0 = 1$~kg, $k_0 = 1$~N~m$^{-1}$ and $c_0 = 0.02$~kg~s$^{-1}$, giving rise to 1\% modal damping. The mass ratio between the absorber and the host structure is 5\%. Fig.~\ref{fig:sdofDampedInit} shows that the peaks of the compliance are unbalanced; the initial tuning is thus to be improved.

            \subsection{Dynamics of a structure with multiple tuned mass dampers}
	    \label{ssec:dynamics}	
	    
		\subsubsection{Dynamics of the uncontrolled structure}

			The dynamics of the discretized host structure is governed by the set of $N$ coupled linear second-order ordinary differential equations (ODEs)
			\begin{equation}
				\mathbf{M}_0\ddot{\mathbf{x}} + \mathbf{C}_0\dot{\mathbf{x}} + \mathbf{K}_0\mathbf{x} = \mathbf{f},
				\label{eq:hostEOMTime}
			\end{equation}
			where $\mathbf{M}_0$, $\mathbf{C}_0$ and $\mathbf{K}_0$ are the structural mass, damping and stiffness matrices, respectively, $\mathbf{x}$ is the vector of generalized coordinates and $\mathbf{f}$ is the vector of conjugated generalized forces. Subscript $0$ refers to the host structure, and an upper dot ($\dot{\ }$) denotes time derivation. Assuming harmonic forcing at angular frequency $\omega$, the Fourier transform of Eq.~\eqref{eq:hostEOMTime} is given by
			\begin{equation}
				\left(-\omega^2\mathbf{M}_{0} + j\omega \mathbf{C}_0 + \mathbf{K}_0\right)\mathbf{X}(\omega) = \mathbf{H}_0(\omega)\mathbf{X}(\omega) = \mathbf{F}(\omega).
				\label{eq:hostEOMFreq}
			\end{equation}
			$j$ is the unit imaginary number (i.e., $j^2 =-1$), $\mathbf{X}$ and $\mathbf{F}$ are the Fourier transforms of $\mathbf{x}$ and $\mathbf{f}$, respectively, and $\mathbf{H}_0$ is the dynamic stiffness matrix. Because the host structure is in general lightly damped, the assumption of proportional damping is made. Using a modal expansion of the displacement and projecting the equations of motion onto the modal basis, the following inverse relation can be derived~\cite{geradin2014mechanical}:
			\begin{equation}
				\mathbf{X}(\omega) = \mathbf{H}_0^{-1}(\omega)\mathbf{F} = \mathbf{\Phi}_0\left(-\omega^2\mathbf{I} + 2j\omega \mathbf{Z}_0\mathbf{\Omega}_0 + \mathbf{\Omega}_0^2\right)^{-1} \mathbf{\Phi}^T_0\mathbf{F}(\omega),
				\label{eq:hostModalExpansion}
			\end{equation}
			where $\mathbf{\Phi}_0$ is the matrix of the mass-normalized mode shapes, $\mathbf{I}$ is the identity matrix, $\mathbf{\Omega}_0$ is a diagonal matrix containing the undamped resonance frequencies of the structure and $\mathbf{Z}_0$ is a diagonal matrix containing the associated modal damping coefficients. The superscript $T$ denotes matrix transposition. In the remainder of this article, the modal expansion given by Eq.~\eqref{eq:hostModalExpansion} is assumed to be known, so that $\mathbf{H}_0^{-1}$ is known as well. 
			
			\subsubsection{Dynamics of the controlled structure}

			The dynamic equations of the controlled structure with $N_a$ absorbers are given by the set of $N+N_a$ ODEs
			\begin{equation}
				\begin{array}{c} \left[ \begin{array}{cc} \mathbf{M}_0 & \mathbf{0} \\ \mathbf{0} & \mathbf{M}_a \end{array} \right] \left[\begin{array}{c} \ddot{\mathbf{x}} \\ \ddot{\mathbf{x}}_a \end{array} \right] + \left[ \begin{array}{cc} \mathbf{C}_0 + \mathbf{B}\mathbf{C}_a \mathbf{B}^T & -\mathbf{B}\mathbf{C}_a \\ -\mathbf{C}_a \mathbf{B}^T & \mathbf{C}_a \end{array} \right] \left[\begin{array}{c} \dot{\mathbf{x}} \\ \dot{\mathbf{x}}_a \end{array} \right]  + \\
				\left[ \begin{array}{cc} \mathbf{K}_0  + \mathbf{B}\mathbf{K}_a \mathbf{B}^T & -\mathbf{B}\mathbf{K}_a \\ -\mathbf{K}_a \mathbf{B}^T & \mathbf{K}_a \end{array} \right] \left[\begin{array}{c} \mathbf{x} \\ \mathbf{x}_a \end{array} \right]  = \left[\begin{array}{c} \mathbf{f} \\ \mathbf{0} \end{array} \right],
				\end{array}
				\label{eq:EOMControlledStructure}
			\end{equation}
			where $\mathbf{x}_a$ is the vector of generalized coordinates associated with the absorbers, $\mathbf{B}$ is a localization matrix collecting every localization vector $\mathbf{b}_n$ associated with the $n^{th}$ absorber
			\begin{equation}
				\mathbf{B} = \left[\mathbf{b}_1\ ,\ \cdots\ ,\ \mathbf{b}_{N_a}\right],
			\end{equation}
			and $\mathbf{M}_a$, $\mathbf{C}_a$ and $\mathbf{K}_a$ are diagonal matrices containing TMD parameters
			\begin{equation}
				\mathbf{M}_a = \left[\begin{array}{ccc} m_{a,1} & & \\ & \ddots & \\ & & m_{a,N_a} \end{array}\right], \qquad \mathbf{C}_a = \left[\begin{array}{ccc} c_{a,1} & & \\ & \ddots & \\ & & c_{a,N_a} \end{array}\right], \qquad \mathbf{K}_a = \left[\begin{array}{ccc} k_{a,1} & & \\ & \ddots & \\ & & k_{a,N_a} \end{array}\right],
				\label{eq:absorbersMatrix}
			\end{equation}
			in which $m_{a,n}$, $c_{a,n}$ and $k_{a,n}$ are the mass, damping and stiffness of the $n^{th}$ absorber, respectively.
			
			Expressing the equations of motion in the frequency domain, it is possible to derive the compliance matrix in a manner similar to Eqs.~\eqref{eq:hostEOMFreq}--\eqref{eq:hostModalExpansion}. The burden associated with computing the compliance may however be alleviated by the  Sherman-Morrison-Woodbury (SMW) formula~\cite{woodbury1950inverting} 
			\begin{equation}
				\left(\mathbf{A} + \mathbf{U}\mathbf{Q}\mathbf{V}\right)^{-1} = \mathbf{A}^{-1} + \mathbf{A}^{-1}\mathbf{U}\left(\mathbf{Q}^{-1} - \mathbf{V}\mathbf{A}^{-1}\mathbf{U}\right)^{-1}\mathbf{V}\mathbf{A}^{-1}.
				\label{eq:SMW}
			\end{equation}
			for invertible matrices $\mathbf{A}$ and $\mathbf{Q}$. The principles of this alleviation were proposed in previous works. Ozer and Royston~\cite{Ozer2005} used the Sherman-Morrison formula~\cite{Sherman1950} as a simplifying numerical tool to tune the parameters of one absorber. A generalization to multiple lumped elements based on the SMW formula was later proposed by Cha~\cite{Cha2007}, but no attempt was made to use the formula to tune the absorbers. 
			
			The Fourier transform of the second line of Eq.~\eqref{eq:EOMControlledStructure} yields 
			\begin{equation}
				\mathbf{X}_a (\omega) = \left(-\omega^2 \mathbf{M}_a + j\omega \mathbf{C}_a + \mathbf{K}_a\right)^{-1}\left(j\omega \mathbf{C}_a + \mathbf{K}_a\right)\mathbf{B}^T\mathbf{X}(\omega).
				\label{eq:absorberCondensation}
			\end{equation}
			Inserting this relation back into the Fourier transform of the first line of Eq.~\eqref{eq:EOMControlledStructure}, one gets
			\begin{equation}
				\left(\mathbf{H}_0(\omega) + \mathbf{B}\left\{j\omega \mathbf{C}_a + \mathbf{K}_a - \left(j\omega \mathbf{C}_a + \mathbf{K}_a\right)\left(-\omega^2 \mathbf{M}_a + j\omega \mathbf{C}_a + \mathbf{K}_a\right)^{-1}\left(j\omega \mathbf{C}_a + \mathbf{K}_a\right)\right\}\mathbf{B}^T\right)\mathbf{X}(\omega) = \mathbf{F}(\omega).
			\end{equation}
			Carrying out simplifications on the diagonal matrices related to the absorbers, the dynamic stiffness matrix of the controlled structure $\mathbf{H}_c$ can be expressed as
			\begin{equation}
				\mathbf{H}_c(\omega) = \mathbf{H}_0(\omega) + \mathbf{B}\mathbf{H}_A(\omega)\mathbf{B}^T.
				\label{eq:HCondensed}
			\end{equation}
			where $\mathbf{H}_A$ is a diagonal matrix given by
			\begin{equation}
				\mathbf{H}_A(\omega) = -\omega^2\left(j\omega \mathbf{C}_a + \mathbf{K}_a\right)\left(-\omega^2 \mathbf{M}_a + j\omega \mathbf{C}_a + \mathbf{K}_a\right)^{-1} \mathbf{M}_a 		
			\end{equation}
			Thus, the dynamic stiffness matrix of the controlled structure is equal to the sum of the dynamic stiffness matrix of the host structure and a rank-$N_a$ update representing the feedback action of the absorbers on the host structure. Consequently, the SMW formula can be used to compute the compliance of the controlled structure as
			\begin{equation}
				\mathbf{H}_c^{-1}(\omega) = \mathbf{H}_0^{-1}(\omega) - \mathbf{H}_0^{-1}(\omega)\mathbf{B} \left(\mathbf{H}_A^{-1}(\omega)+\mathbf{B}^T\mathbf{H}_0^{-1}(\omega)\mathbf{B}\right)^{-1} \mathbf{B}^T\mathbf{H}_0^{-1}(\omega)
				\label{eq:SMWFRF}
			\end{equation} Eq.~\eqref{eq:SMWFRF} can be subject to singularity issues in three cases. The first one is $\omega=0$, because $\mathbf{H}_A(0)$ is a zero matrix. In that case, Eq.~\eqref{eq:HCondensed} simply indicates that $\mathbf{H}^{-1}_c(0) = \mathbf{H}_0^{-1}(0)$. The second case occurs if any $m_{a,n}$ is zero, or if any pair $(c_{a,n},k_{a,n})$ is zero. These cases correspond to an absence of absorber or to an unattached absorber mass, which is irrelevant in the design problem. Finally, the SMW formula requires $\mathbf{H}_0$ to be non-singular. This condition might not be met at the resonance frequencies of an undamped host structure. In this case, a small amount of damping may be added to resolve this numerical issue while still representing faithfully the dynamics of the host.			
			
			Finally, the compliance at a given displacement located by a vector $\mathbf{w}_u$, is 
				\begin{equation}
					h(\omega) = \mathbf{w}_u^T\mathbf{H}_c^{-1}(\omega)\mathbf{w}_f,
				    \label{eq:compliance}
				\end{equation}				    		
    			where $\mathbf{w}_f$ is a vector describing the spatial distribution of the forcing vector $\mathbf{F}$, and $\mathbf{H}_c^{-1}$ is evaluated using Eq.~\eqref{eq:SMWFRF}. Alternatively, the accelerance (i.e., the transfer function between an acceleration of interest and the external forcing amplitude) may be considered simply as
    			\begin{equation}
    			    h_a(\omega) = -\omega^2h(\omega)
    			\end{equation}
			
		
    		 
    		\subsection{Peak finding}
    		\label{ssec:peaks}

				The resonance frequencies occur at the maximum of the compliance, i.e.,  
				\begin{equation}
					\omega_s = \arg \min_{\omega \in \mathbb{R}} s_1(\omega) = \arg \min_{\omega \in \mathbb{R}}\left(-|h(\omega)|^2\right),
					\label{eq:s1}
				\end{equation}				    		
    			where the square modulus of the complex compliance is used to make the function $s_1$ smooth with respect to $\omega$. A necessary condition to satisfy this relation is
				\begin{equation}
					\omega_s : \left. \frac{\partial s_1(\omega)}{\partial \omega}\right|_{\omega=\omega_s} = - 2\left.\left(\frac{\partial h^*(\omega)}{\partial \omega}h(\omega) + h^*(\omega)\frac{\partial h(\omega)}{\partial \omega}\right)\right|_{\omega=\omega_s} = 0,
					\label{eq:s1derivative}
				\end{equation}    		
    			where the superscript $*$ denotes a complex conjugation. This equation can be solved numerically starting from an initial guess (e.g. Eq.~\eqref{eq:peaksInit}) using either root-finding algorithms (paying attention to the fact that a root might not correspond to a maximum of the compliance) or linesearch algorithms~\cite{nocedal2006numerical}. This procedure yields a set of frequencies $\omega_{i}$ and associated amplitudes noted $|h(j\omega_{i})| = |h|_i$ with $i=1,...,2N_a$. 
    			
    			A similar procedure for the accelerance consists in finding the roots of
    			\begin{equation}
					\omega_s : \left. \frac{\partial s_2(\omega)}{\partial \omega}\right|_{\omega=\omega_s} = -4\omega^3\left| h(j\omega)\right|^2 - 2\omega^4\left.\left(\frac{\partial h^*(\omega)}{\partial \omega}h(\omega) + h^*(\omega)\frac{\partial h(\omega)}{\partial \omega}\right)\right|_{\omega=\omega_s} = 0,
					\label{eq:s2derivative}
				\end{equation}   
    			

    		\subsection{$p$-norm optimization}
    		\label{ssec:pnorm}
    		
    		The goal of the optimization algorithm is to find the optimal mass, damping and stiffness of the absorbers through the nonlinear programming problem
            \begin{equation}
					\begin{aligned}
					& \underset{\bm{\xi}}{\text{minimize}}
					& & f_p(\xi) \\
					& \text{subject to}
					& & \sum_{n=1}^{N_a} m_{a,n} - m_{\text{Max}}\leq 0
					\end{aligned},
					\label{eq:optimizationProblem}
				\end{equation}				
			where $\xi$ is the vector containing the absorber parameters, and $f_p$ is the $p$-norm of the vector containing the squared amplitudes of the controlled resonance peaks 
				\begin{equation}
					f_p = \chi\sqrt[p]{\sum_{i=1}^{2N_a} \left(\frac{1}{\chi} |h|^2_{i}\right)^p},
					\label{eq:pnorm}
				\end{equation}	 
    		$\chi$ is a normalizing factor, which does not affect the norm and avoids bad numerical conditioning for large $p$. A typical choice for $\chi$ is 
				\begin{equation}
					\chi = \max_{i\in [1,2N_a]} |h|^2_i
					\label{eq:pnormScaling}
				\end{equation}
			For practical reasons, the total mass of the absorbers should not exceed a maximum $m_{\text{Max}}$, which is translated by the addition of a linear inequality constraint in problem (\ref{eq:optimizationProblem}). It was generally observed in the literature (e.g.~\cite{nishihara2002closed,Leung2008}), and by the authors as well, that the mass constraint is usually active in the optimum design. The convergence of the algorithm may thus be accelerated when an equality constraint is imposed.
    		
			The gradients of the $p$-norm are computed analytically in the proposed algorithm. From Eq.~\eqref{eq:pnorm}, the $l^{th}$ element of the gradient of the $p$-norm with respect to the absorbers parameters is given by
				\begin{equation}
					\dfrac{\partial f_p}{\partial \xi_l} =  \left(\sum_{i=1}^{2N_a} \left(\frac{1}{\chi} |h|^2_{i}\right)^{p-1} \left(\dfrac{\partial h_i^*}{\partial \xi_l}h_i + h_i^*\dfrac{\partial h_i}{\partial \xi_l}\right) \right)\left(\sum_{i=1}^{2N_a} \left(\frac{1}{\chi} |h|^2_{i}\right)^p\right)^{\frac{1}{p}-1}.
					\label{eq:fpGradient}
				\end{equation}
    		The derivatives of the compliance with respect to $\xi_l$ are computed thanks to Eq.~\eqref{eq:SMWFRF} and~\eqref{eq:compliance} as
    			\begin{equation}
    				\dfrac{\partial h_i}{\partial \xi_l} = \mathbf{w}_u^T\mathbf{G}(\omega_i)\dfrac{\partial \mathbf{H}_A^{-1}(\omega_i)}{\partial \xi_l}\mathbf{G}^T(\omega_i)\mathbf{w}_f
    			\label{eq:SMWGradient}
    			\end{equation}
    		where $\omega_i$ are the solutions of Eq.~\eqref{eq:s1derivative} or Eq.~\eqref{eq:s2derivative}, and
				\begin{equation}
				 \mathbf{G}(\omega_i) = \mathbf{H}_0^{-1}(\omega_i)\mathbf{B}\left(\mathbf{H}_A^{-1}(\omega_i)+\mathbf{B}^T\mathbf{H}_0^{-1}(\omega_i)\mathbf{B}\right)^{-1}.
				 \label{eq:SMWGradientG}
				\end{equation}				    			
    		Despite the rather complicated structure of Eqs.~\eqref{eq:SMWGradient} and~\eqref{eq:SMWGradientG}, computing the gradient of the cost function is not cumbersome for two reasons. First, each element in Eq.~\eqref{eq:SMWGradientG} is known from the computation of $|h|_i$. Second, the derivative of $\mathbf{H}_A^{-1}$ with respect to $\xi_l$ can be computed analytically and contains only one non-zero entry. Assuming $\xi_l$ is a parameter associated with the $n^{th}$ absorber, the corresponding entry is given by
    			\begin{equation}
    				\def\arraystretch{2}
    				\left(\dfrac{\partial \mathbf{H}_A^{-1}(\omega_i)}{\partial \xi_l}\right)_{n,n} = \left\{\begin{array}{cl} 
    				\displaystyle \dfrac{1}{m_{a,n}^2 \omega_i^2} & \qquad \text{if}\ \xi_l = m_{a,n} \\
    				\displaystyle -\dfrac{j\omega_i}{\left(k_{a,n}+j\omega_i c_{a,n}\right)^2}& \qquad \text{if}\ \xi_l = c_{a,n} \\ 
    				\displaystyle -\dfrac{1}{\left(k_{a,n}+j\omega_i c_{a,n}\right)^2}& \qquad \text{if}\ \xi_l = k_{a,n}  
    				\end{array}\right.
    			\label{eq:HAGradient}
    			\end{equation}
    		Hence, the gradients of the cost function are obtained analytically by plugging Eq.~\eqref{eq:HAGradient} into Eq.~\eqref{eq:SMWGradient} and then into Eq.~\eqref{eq:fpGradient}. The gradients of a cost function based on the accelerance are obtained through multiplication by $\omega_i^4$.
    			
\begin{figure}[!ht]
				\centering
				\includegraphics[width=0.6\textwidth]{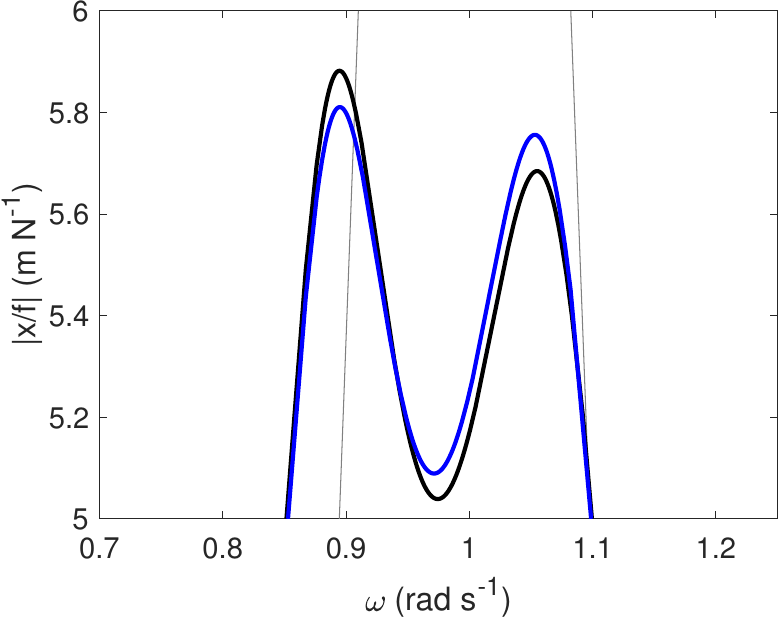}
				\caption{Compliance of the damped single-degree-of-freedom system: uncontrolled structure (\textcolor{gray}{---}), initial tuning (\textbf{---}) and $p$-norm optimization with $p=1$ (\textcolor{blue}{\textbf{---}}).}
				\label{fig:sdofDamped}
			\end{figure}

    		The result of the $1$-norm optimization of the initial tuning in Fig. \ref{fig:sdofDampedInit} is displayed in Fig.~\ref{fig:sdofDamped}. The algorithm has thus reduced the initial mistuning.	
    			
    		\subsection{Norm-homotopy optimization procedure}\label{norm}
    		
				Once the optimization has converged for a given value of $p$, $p$ is then increased in order to penalize high-amplitude peaks more strongly and approach the $H_\infty$ optimum. A heuristic scheme for $p$ given by the double exponential progression
				\begin{equation}
					p = 2^{2^k}, \qquad k\in \mathbb{N}.
					\label{eq:pFunctionk}
				\end{equation}
				is considered. The value of $k$ starts from zero and is incremented by one after convergence. This norm-homotopy algorithm may be terminated when no significant change is observed in the absorbers parameters and/or in the value of the p-norm.
				
\begin{figure}[!ht]
				\centering
				\includegraphics[width=0.6\textwidth]{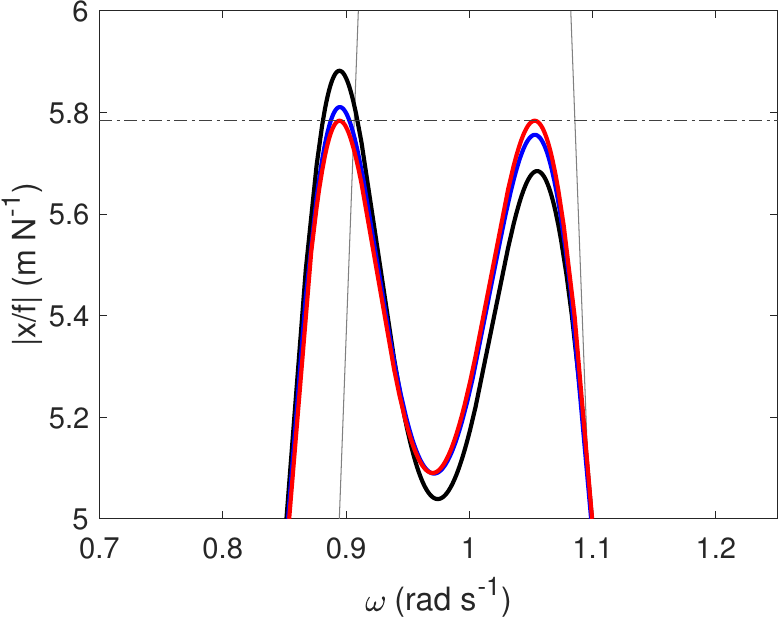}
				\caption{TMD on a damped single-degree-of-freedom system: uncontrolled structure (\textcolor{gray}{---}), initial tuning (\textbf{---}), $p$-norm optimization with $p=1$ (\textcolor{blue}{\textbf{---}}) and norm-homotopy optimization (\textcolor{red}{\textbf{---}}).}
				\label{fig:sdofDampednh}
			\end{figure}

			The end result of the norm-homotopy procedure applied to the single-degree-of-freedom system is shown in Fig.~\ref{fig:sdofDampednh}. The optimization was stopped when $p$ was equal to 65536. The peaks are now equal. 
			
			Two remarks are outlined to punctuate this section. First, the compliance and accelerance were used throughout this work in the cost function of the optimization algorithm, but other transfer functions such as mobility could easily be considered as well. Second, only one input-output pair was taken into account (through the vectors $\mathbf{w}_u$ and $\mathbf{w}_f$) for this transfer function for simplicity; more elaborate cost functions than that given in Eq.~\eqref{eq:pnorm} could be conceived to generalize the approach to multiple inputs and/or outputs (see e.g.~\cite{Salcedo-Sanz2017}).

    	\section{Examples}
    	\label{sec:examples}
			Two examples serve to demonstrate the algorithm in this section. The maximum allowable mass for the absorbers is set to 5\% of that of the host structure. The optimization problem \eqref{eq:optimizationProblem} is solved in MATLAB thanks to the \texttt{fmincon} routine. Each time the cost function is called, the peak-finding algorithm solves Eqs.~\eqref{eq:s1}-\eqref{eq:s1derivative} using a linesearch approach, as described in~\cite{nocedal2006numerical}.    	
    	
    		\subsection{A two-degree-of-freedom host structure}
    		\label{ssec:2dofs}
    		
\begin{figure}[!ht]
						\centering
						\includegraphics[width=0.6\textwidth]{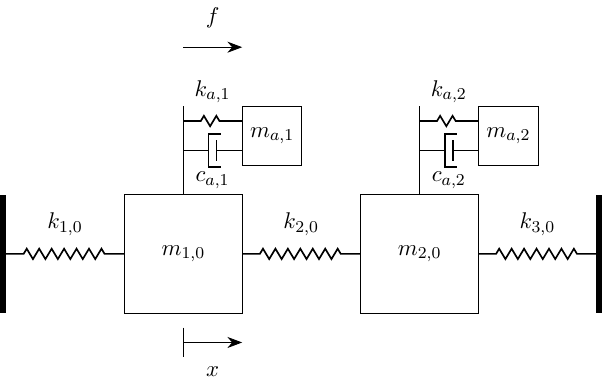}
						\caption{Schematic representation of the two-degree-of-freedom structure with the two attached absorbers.}
						\label{fig:2dofSch}
					\end{figure}

				The two-degree-of-freedom structure with the two attached absorbers is depicted in Fig.~\ref{fig:2dofSch}. The parameters of the host system are listed in Table~\ref{tab:2dofChar}. The first absorber, labelled "1" in Fig.~\ref{fig:2dofSch}, targets the first mode whereas the second absorber targets the second mode.    		

				\begin{table}[!ht]
					\centering
					\begin{tabular}{rccccc}
						\hline
							Parameter & $m_{1,0}$ (kg) & $m_{2,0}$ (kg) & $k_{1,0}$ (N m$^{-1}$) & $k_{2,0}$ (N m$^{-1}$) & $k_{3,0}$ (N m$^{-1}$) \\
						\hline
							Value & 1 & 1 & 1 & 1 & 1 \\
						\hline
					\end{tabular}
					\caption{Parameters of the two-degree-of-freedom host structure.}
					\label{tab:2dofChar}
				\end{table}	
							    		
\begin{figure}[!ht]
						\centering
						\includegraphics[width=0.6\textwidth]{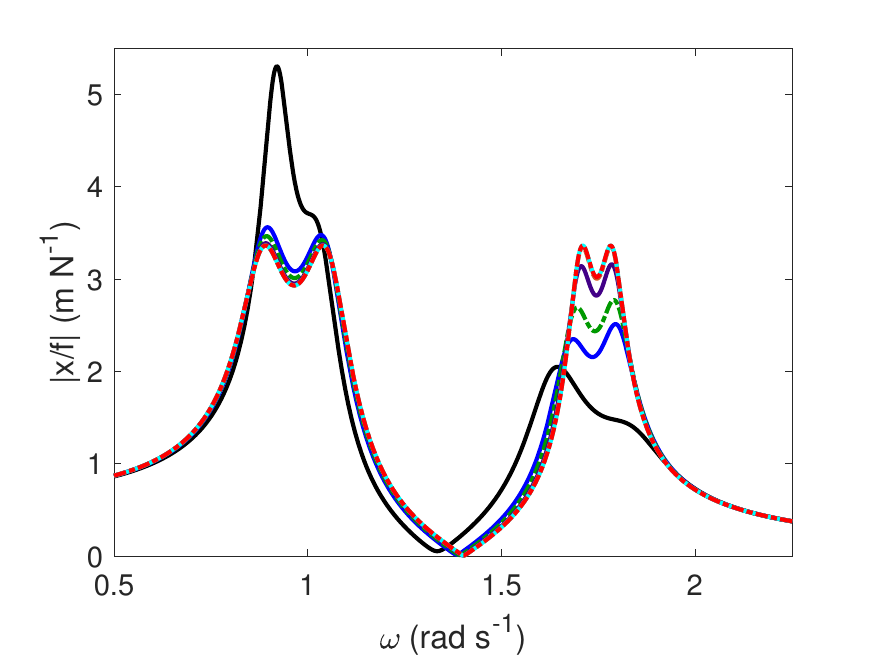}
						\caption{MTMD on a two-degree-of-freedom system: initial tuning (\textbf{---}), solution for $k=0$ (\textcolor{blue}{\textbf{---}}), solution for $k=1$ (\textcolor{darkGreen1}{\textbf{-$\cdot$-}}), solution for $k=2$ (\textcolor{violet1}{\textbf{---}}), solution for $k=3$ (\textcolor{orange1}{\textbf{-$\cdot$-}}), solution for $k=4$ (\textcolor{cyan}{\textbf{---}}), and norm-homotopy optimal solution (\textcolor{red}{\textbf{-$\cdot$-}}).}
						\label{fig:2dofs}
					\end{figure}

				Fig.~\ref{fig:2dofs} displays the compliance of the host structure for different values of $k$. The initial tuning is clearly unsatisfactory. After the optimization for $k=0$, almost equal peaks around the two resonances are retrieved, but the amplitude of the peaks around the first mode is still significantly larger than that around the second mode. Increasing $k$ up to 4 eventually leads to the so-called all-equal-peak design. The optimal parameters of the absorbers are listed in Table~\ref{tab:2dofAbsChar}.
    		
    				\begin{table}[!ht]
					\centering
					\begin{tabular}{cccc}
						\hline
						Parameter & $m_a$ (kg) & $c_a$ (N s m$^{-1}$) & $k_a$ (N m$^{-1}$) \\ 
						\hline
						Absorber 1 & 0.94$m_{\text{Max}}$ & 0.0237  & 0.0840 \\
						Absorber 2 & 0.06$m_{\text{Max}}$ & 0.0007 &  0.0183 \\
						\hline
					\end{tabular}
					\caption{Optimal parameters of the absorbers in Fig.~\ref{fig:2dofSch} for compliance optimization.}
					\label{tab:2dofAbsChar}
				\end{table}
				
				\subsubsection{Accelerance optimization}
				
\begin{figure}[!ht]
						\centering
						\includegraphics[width=0.6\textwidth]{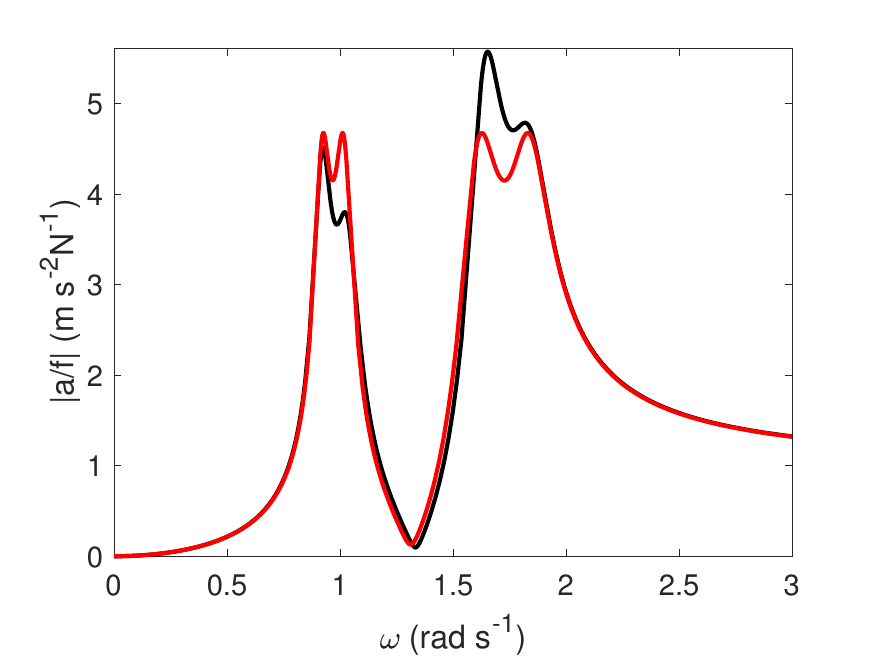}
						\caption{MTMD on a two-degree-of-freedom system: initial tuning (\textbf{---}), and norm-homotopy optimal solution (\textcolor{red}{\textbf{-$\cdot$-}}).}
						\label{fig:2dofsAcc}
					\end{figure}

				Figure~\ref{fig:2dofsAcc} shows the results of a norm-homotopy optimization on the accelerance, and Table~\ref{tab:2dofAbsCharAcc} gives the corresponding optimal characteristics. Again, peaks of equal amplitude are observed, which illustrates the versatility of the approach. 
				
					\begin{table}[!ht]
					\centering
					\begin{tabular}{cccc}
						\hline
						Parameter & $m_a$ (kg) & $c_a$ (N s m$^{-1}$) & $k_a$ (N m$^{-1}$) \\ 
						\hline
						Absorber 1 & 0.37$m_{\text{Max}}$ & 0.0057  & 0.0342 \\
						Absorber 2 & 0.63$m_{\text{Max}}$ & 0.0213 &  0.1720 \\
						\hline
					\end{tabular}
					\caption{Optimal parameters of the absorbers in Fig.~\ref{fig:2dofSch} for accelerance optimization.}
					\label{tab:2dofAbsCharAcc}
				\end{table}

    		\subsection{A simply-supported aluminum plate}
			\label{ssec:plate}
			    		
    			The second example is a homogeneous, isotropic, simply-supported rectangular plate that features closely-spaced resonances. According to Kirchhoff-Love theory, the mode shapes and eigenfrequencies of a plate of length $a$, width $b$ and thickness $h$ are given by
    			\begin{equation}
    				\begin{array}{rl}
    				\phi_{mn}(x,y) & = \dfrac{2}{\sqrt{M}}\sin \left(\dfrac{m\pi x}{a}\right)\sin \left(\dfrac{n\pi y}{b}\right)\\
    				\omega_{mn} &= \sqrt{\dfrac{D}{\rho h}}\left[\left(\dfrac{m\pi}{a}\right)^2 + \left(\dfrac{n\pi}{b}\right)^2  \right]
    				\end{array},
    			\end{equation}
			    respectively~\cite{geradin2014mechanical}. $\rho$ is the density of the plate, $M=\rho a b h$ is the mass of the plate, and $D = Eh^3/(12(1-\nu^2))$ is the plate bending stiffness, where $E$ and $\nu$ are Young's modulus Poisson's ratio, respectively. The plate parameters are given in Table~\ref{tab:plateChar}.
				
				\begin{table}[!ht]
					\centering
					\begin{tabular}{cc}
						\hline						
						\textbf{Characteristic} & \textbf{Value}\\						
						\hline
						Length $a$ & 1 m \\ 
						Width $b$ & 0.7 m\\
						Thickness $h$ & 1 mm \\
						Young modulus $E$ & 68 GPa \\
						Poisson ratio $\nu$ & 0.36 \\
						Density $\rho$ & 2700 kg m$^{-3}$\\
						\hline
					\end{tabular}
					\caption{Parameters of the simply-supported aluminum plate.}
					\label{tab:plateChar}
				\end{table}
				
				To discretize the model, the mode shapes are sampled spatially at locations $(\mathbf{x}_s,\mathbf{y}_s)$, where
				    \begin{equation}
				        \begin{array}{ccc}
				             \mathbf{x}_s & = \left[x_u,x_{a,1},x_{a,2},x_{a,3},x_{a,4},x_f\right]^T &= \left[0.25a, 0.5a, 0.15a, 0.4a ,0.25a , 0.75a\right]^T  \\
				             \mathbf{y}_s & = \left[y_u,y_{a,1},y_{a,2},y_{a,3},y_{a,4},y_f\right]^T &= \left[0.25b, 0.5b, 0.4b, 0.15b, 0.75b, 0.75b\right]^T 
				        \end{array},
				    \end{equation}
				and only a finite number of modes is retained, up to $m=M_{\text{Max}}=10$ and $n=N_{\text{Max}}=10$. The mode shape matrix of the host structure such that  
				\begin{equation}
					\mathbf{\Phi}_0 = \left[\phi_{11}(\mathbf{x}_s,\mathbf{y}_s),\phi_{12}(\mathbf{x}_s,\mathbf{y}_s), \cdots ,\phi_{M_{\text{Max}}N_{\text{Max}}}(\mathbf{x}_s,\mathbf{y}_s)\right],
				\end{equation}
				\begin{equation}
					\mathbf{\Omega}_0^2 = \mbox{diag}\left(\omega_{11}^2,...,\omega_{M_{\text{Max}}N_{\text{Max}}}^2 \right).
				\end{equation}
				
\begin{figure}[!ht]
					\centering
					\includegraphics[width=0.6\textwidth]{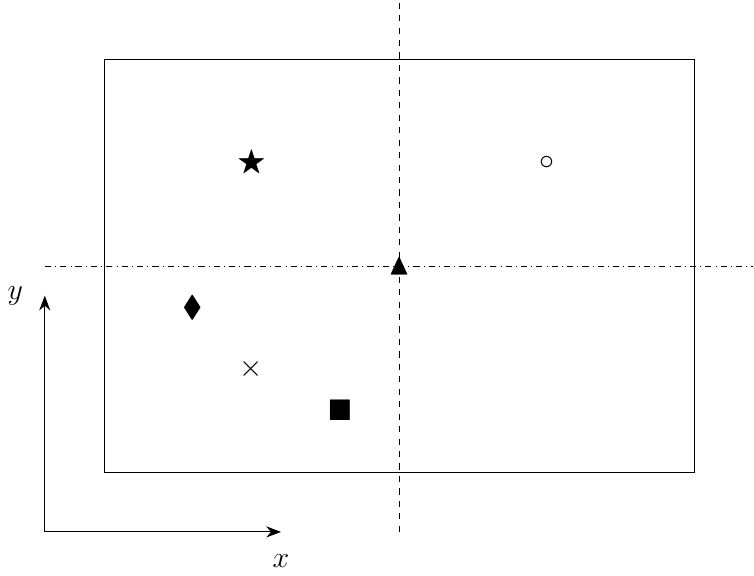}
					\caption{Geometry of the plate: point force location ($\times$), displacement measurement location ($\circ$), first TMD location ($\blacktriangle$), second TMD location ($\blacklozenge$), third TMD location ($\blacksquare$), fourth TMD location ($\bigstar$), nodal line of modes for which $m=2$ (- -) and nodal line of modes for which $n=2$ ($-\cdot -$).}
					\label{fig:plateGeom}
				\end{figure}

				The plate is loaded by a harmonic point force located at $(x_f,y_f)$. Three/four absorbers are considered to mitigate the first three/four resonances, respectively. Fig.~\ref{fig:plateGeom} depicts the geometrical configuration. 
				
				To have a numerically well-conditioned problem, the compliance measured  at coordinates $(x_u,y_u)$ is normalized with the static displacement $x_{\text{st}}$ of the structure 
				\begin{equation}
					x_{\text{st}}(x_u,y_u) = \sum_{m=1}^{M_{\text{Max}}}\sum_{n=1}^{N_{\text{Max}}} \dfrac{\phi_{mn}(x_u,y_u)\phi_{mn}(x_f,y_f)}{\omega_{mn}^2}f.
				\end{equation}

			    The results of a norm-homotopy optimization are compared to those of a direct $H_\infty$ norm optimization. The $H_\infty$ norm was computed using a standard method~\cite{Bruinsma1990} and without using the SMW formula.
				
				\subsubsection{Vibration mitigation with three absorbers}			
			
				The first three modes, i.e., $(m,n)$ = $(1,1)$, $(2,1)$ and $(1,2)$, are first targeted. The first absorber is placed at the antinode of the first mode, whereas the second and third absorbers are slightly offset in order to influence higher-frequency modes which are likely to have a nodal line passing through the antinodes of the second and third modes.
				
\begin{figure}[!ht]
					\centering
					\includegraphics[width=0.8\textwidth]{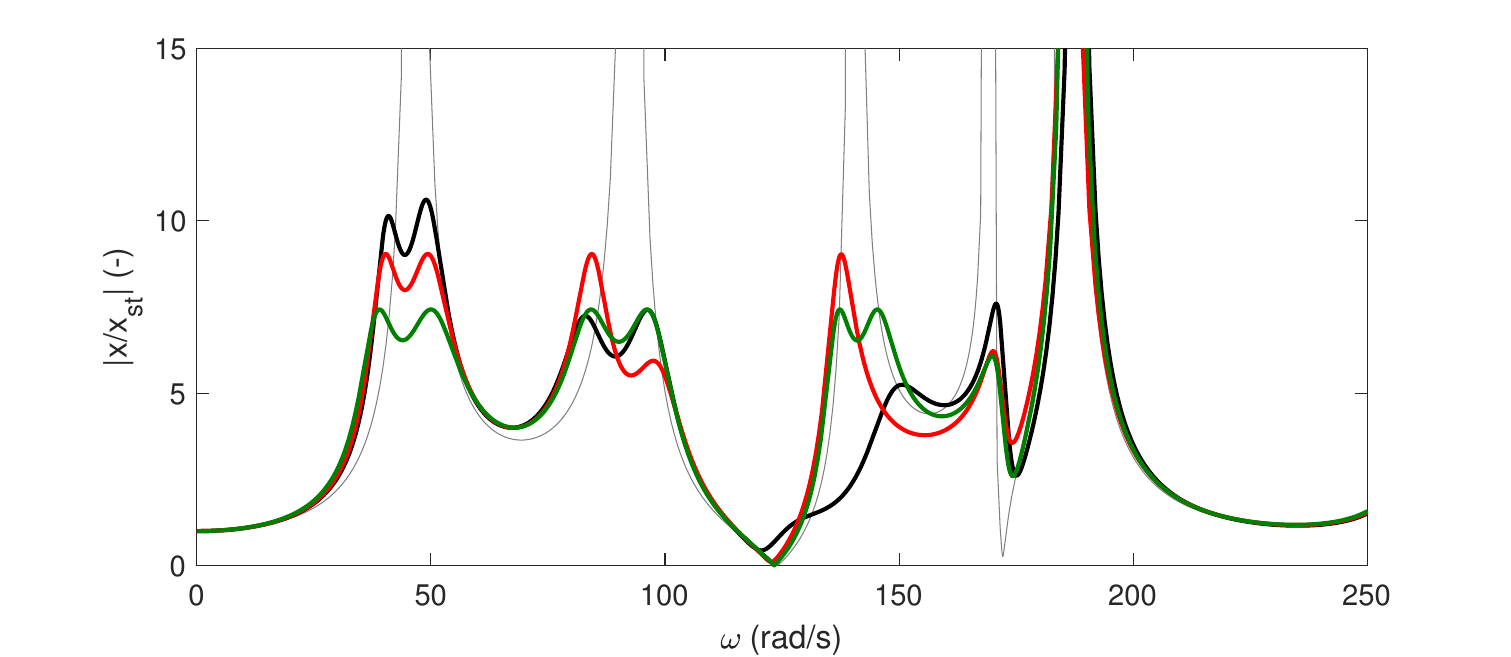}
					\caption{Compliance of the plate with three absorbers targeting modes $(1,1)$, $(2,1)$ and $(1,2)$:  uncontrolled structure (\textcolor{gray}{---}), initial tuning (\textbf{---}) and optimized tuning (\textcolor{red}{\textbf{---}}: direct $H_\infty$ optimization, \textcolor{darkGreen1}{\textbf{---}}: norm-homotopy optimization).}
					\label{fig:plate3abs}
				\end{figure}

				The result of the $H_\infty$ optimization (limiting the range of frequencies up to $\omega=150$ rad s$^{-1}$) and the norm-homotopy optimization are presented in Fig.~\ref{fig:plate3abs}, and Table~\ref{tab:plate3abs} lists the parameters of the resulting absorbers. Although it shows improvement with respect to the initial tuning, the direct $H_\infty$ norm optimization stops somewhat prematurely, and only four peaks are equated in amplitude. With the norm-homotopy optimization, the six peaks around the first three resonances all feature the same amplitude, which further validates the proposed methodology. 
				
				\begin{table}[!ht]
					\centering
					\begin{tabular}{cccc}
						\hline
						\textbf{Absorber number}  & $m_a$ (kg) & $c_a$ (N s m$^{-1}$) & $k_a$ (N m$^{-1}$) \\
						\hline
						1 & 0.65$m_{\text{Max}}$& 1.1348 & 110.42\\
						2 & 0.3$m_{\text{Max}}$& 0.6044 & 225.36 \\
						3 & 0.05$m_{\text{Max}}$ & 0.0666 & 94.93\\
						\hline
					\end{tabular}
					\caption{Parameters of the three absorbers.}
					\label{tab:plate3abs}
				\end{table}
				
				\subsubsection{Vibration mitigation with four absorbers}				
				
\begin{figure}[!ht]
					\centering
					\includegraphics[width=0.8\textwidth]{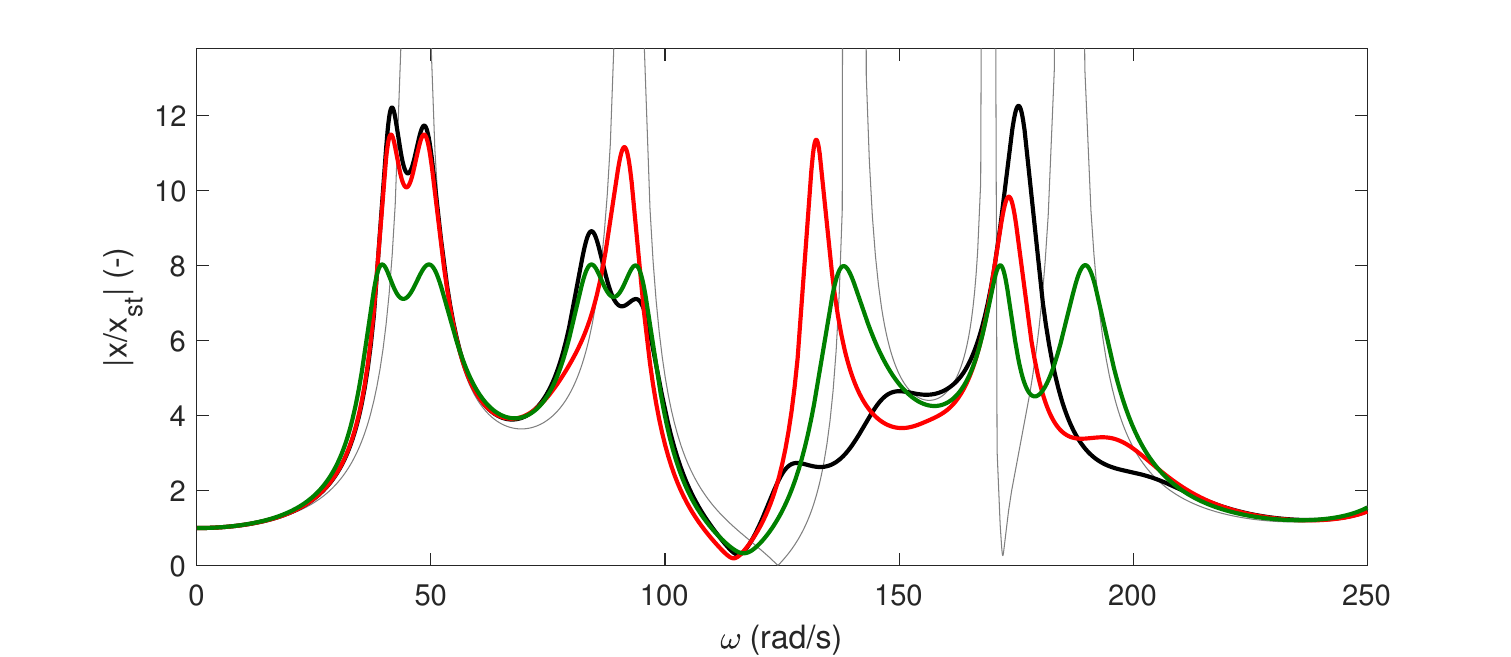}
					\caption{Compliance of the plate with four absorbers targeting modes $(1,1)$, $(2,1)$, $(1,2)$ and $(2,2)$: uncontrolled structure (\textcolor{gray}{---}), initial tuning (\textbf{---}) and optimized tuning (\textcolor{red}{\textbf{---}}: direct $H_\infty$ optimization, \textcolor{darkGreen1}{\textbf{---}}: norm-homotopy optimization).}
					\label{fig:plate4abs}
				\end{figure}

				In the previous example, the fifth mode of the plate remained largely unaffected. To improve the situation, a fourth absorber targeting this mode is placed on the plate. Fig.~\ref{fig:plate4abs} displays the result of the optimization processes (where a maximum frequency of $\omega=250$ rad s$^{-1}$ was considered when computing the $H_\infty$ norm). The direct $H_\infty$ norm optimization only features marginal improvement compared to the initial design. The norm-homotopy optimization obeys the all-equal-peak design and has a lower $H_\infty$ norm, in the considered frequency range. However, we note that modes 3, 4 and 5 do not feature two peaks around their uncontrolled resonance. Looking at the absorber parameters in Table~\ref{tab:plate4abs} reveals that the third absorber is in fact eliminated by the optimization algorithm (zero mass), while the fourth one is controlling modes 3-5 simultaneously. This result probably originates from the fact that plates have closely-spaced frequencies and that the absorber is located relatively far from any nodal line of these modes.

				\begin{table}[!ht]
					\centering
					\begin{tabular}{cccc}
						\hline
						\textbf{Absorber number} & $m_a$ (kg) & $c_a$ (N s m$^{-1}$) & $k_a$ (N m$^{-1}$) \\
						\hline
						1 & 0.56$m_{\text{Max}}$& 0.9074 & 96.0221\\
						2 & 0.22$m_{\text{Max}}$ & 0.3665 & 164.7649\\
						3 & 1$\times 10^{-7} m_{\text{Max}}$& 67.1701 & 442.2204\\
						4 & 0.22$m_{\text{Max}}$& 1.0619 & 444.1343 \\
						\hline
					\end{tabular}
					\caption{Parameters of the four absorbers.}
					\label{tab:plate4abs}
				\end{table}
				
			\subsubsection{Computational cost}	
			    The computational cost of the two optimization approaches cannot be easily compared due to the fact that the cost functions are not identical. However, Table~\ref{tab:plateSpeed} seems to indicate that, while the norm-homotopy approach requires more cost function evaluations, the significant speedup provided by the SMW formula enables a faster execution than a direct $H_\infty$ optimization.
			    
			    \begin{table}[!ht]
					\centering
					\begin{tabular}{ccc}
						\hline
						\textbf{Case} & Average time per cost function evaluation & Total time  \\
						\hline
						3 absorbers, NH & 1 & 1563 \\
						3 absorbers, $H_\infty$ & 17.6 & 2956 \\
						4 absorbers, NH & 1.1 & 1097 \\
						4 absorbers, $H_\infty$ & 18.9 & 1342 \\
						\hline
					\end{tabular}
					\caption{Normalized times (with respect to the average time per cost function evaluation for the norm-homotopy optimization of three absorbers) of the different cases (NH stands for norm-homotopy optimization, $H_\infty$ stands for direct $H_\infty$ optimization).}
					\label{tab:plateSpeed}
				\end{table}
			\subsubsection{Design robustness}
				
\begin{figure}[!ht]
					\centering
					\includegraphics[width=0.6\textwidth]{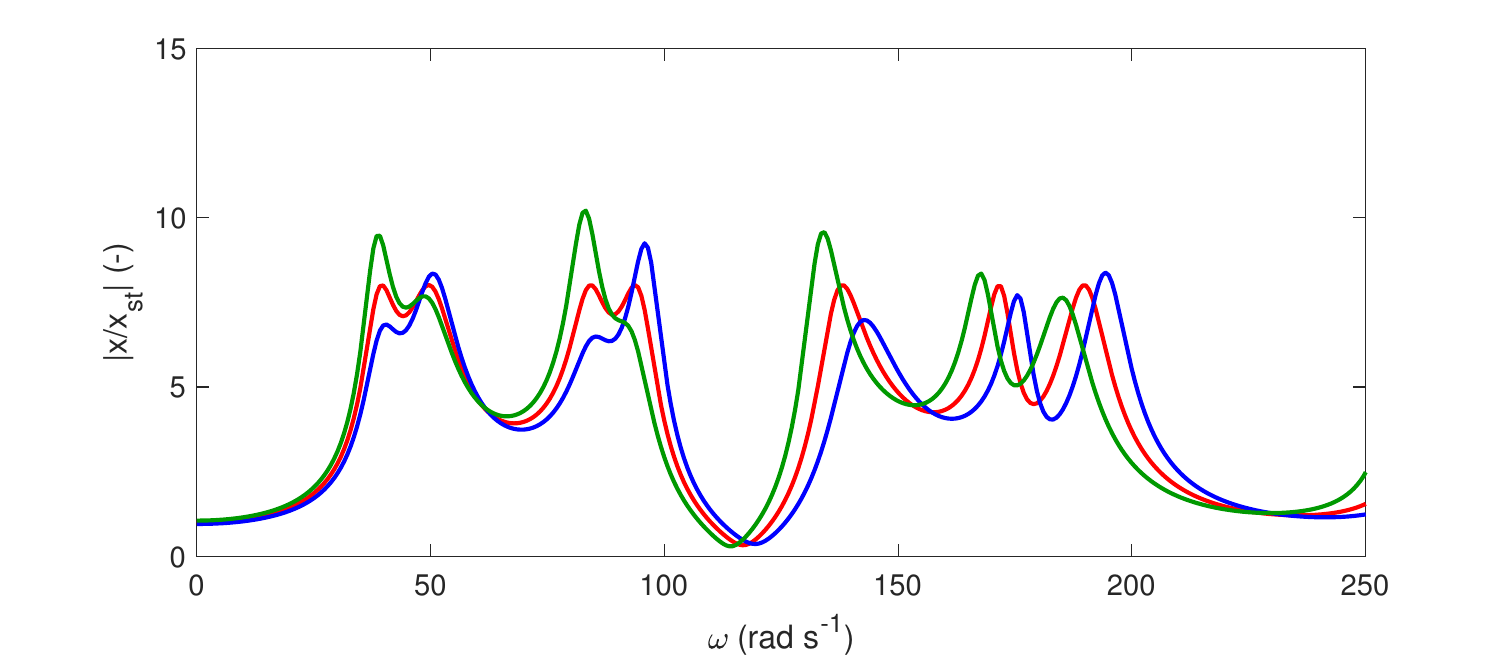}
					\caption{Design robustness: nominal value of $E$ (\textbf{\textcolor{red}{---}}), 5\% increase (\textbf{\textcolor{blue}{---}}) and 5\% decrease  (\textbf{\textcolor{darkGreen1}{---}}).}
					\label{fig:plate4absRobustness}
				\end{figure}

			In real-life applications, the model parameters may be known with limited accuracy. For illustration, variations of $\pm$5\% of Young's modulus are presented herein, while every other parameter is kept constant. The  absorbers parameters of Table~\ref{tab:plate4abs} are used. As depicted in Fig.~\ref{fig:plate4absRobustness}, the absorbers are detuned in a fashion similar to that of the single-degree-of-freedom case, and an increase in the maximum amplification of 15\% (-5\% case) and 27\% (5\% case) can be noticed.
				
    \section{Conclusion}    	
    \label{sec:conclusion}	
		A norm-homotopy numerical optimization algorithm was proposed in this paper to tune multiple TMDs targeting several resonances of a host structure. The algorithm solves a sequence of optimization problems of increasing complexity in which the cost function depends on the $p$-norm of the peak amplitudes of the compliance/accelerance of the controlled structure. As demonstrated by the examples, the outcome of the algorithm is an all-equal-peak design, in which all controlled peaks feature the same amplitude. The algorithm was also found to outperform direct $H_\infty$ optimization and can deal with a variety of discretized structures with moderate computational cost.
		
		Future works may adapt the algorithm to other types of absorbers (e.g., nonlinear absorbers) and may involve the experimental validation of the proposed design approach.
		
		
    \section*{Acknowledgements}
    	
    	Funding: This work was supported by the SPW [WALInnov grant 1610122].
    	
	\appendix
    
    \section{Single TMD on a multiple-degree-of-freedom structure}
    \label{anx:TMD}        
            
            It is considered that only the $n^{th}$ absorber is attached to the host structure. Around its resonance frequency $\omega_r$, it can be assumed that a resonant mode $r$ dominates the structural response. Thus, the following relation approximately holds for $\omega \approx \omega_r$
			\begin{equation}
				\left[\begin{array}{cc} \mathbf{X}(\omega) \\ X_{a,n}(\omega)\end{array}\right] = \left[\begin{array}{cc} \mathbf{\phi}_r & \mathbf{0} \\ 0 & 1\end{array}\right]\left[\begin{array}{cc} \eta_r(\omega) \\ X_{a,n}(\omega)\end{array}\right] = \mathbf{A}_r \left[\begin{array}{cc} \eta_r(\omega) \\ X_{a,n}(\omega)\end{array}\right],
				\label{eq:modalProj}
			\end{equation}
			where $\mathbf{\phi}_r$ is the mass-normalized resonant mode shape, $\eta_r$ is the associated resonant modal coordinate and $X_{a,n}$ is the generalized degree of freedom associated with the $n^{th}$ absorber. Substituting Eq.~\eqref{eq:modalProj} into Eq.~\eqref{eq:EOMControlledStructure} (where only the $n^{th}$ absorber is considered), premultiplying it by $\mathbf{A}_r^T$ and taking into account the modal orthogonality relationships~\cite{geradin2014mechanical}, one gets
			\begin{equation}
			\begin{array}{c}
				\left(-\omega^2 \left[\begin{array}{cc} 1 & 0 \\ 0 & m_{a,n} \end{array}\right] + j\omega \left[\begin{array}{cc} 2\zeta_r \omega_r+\phi_{a,n}^2 c_{a,n} & -\phi_{a,n} c_{a,n} \\ -\phi_{a,n} c_{a,n} & c_{a,n} \end{array}\right] + \right. \\ \left. \left[\begin{array}{cc} \omega_r^2+\phi_{a,n}^2 k_{a,n} & -\phi_{a,n} k_{a,n} \\ -\phi_{a,n} k_{a,n} & k_{a,n} \end{array}\right]\right)				\left[\begin{array}{c} \eta_r(\omega) \\ X_{a,n}(\omega) \end{array}\right] = \left[\begin{array}{c} \mathbf{\phi}_r^T \mathbf{F}(\omega) \\ 0 \end{array}\right],
				\end{array}
				\label{eq:projAbs}
			\end{equation}						
			where $\zeta_r$ is the $r^{th}$ modal damping ratio and $\phi_{a,n} = \mathbf{b}_n^T\mathbf{\phi}_r$ is the mode shape of the host structure at the location where the absorber is to be attached. The base displacement of the $n^{th}$ absorber $U_{n}$ is given by
			\begin{equation}
				U_{n}(\omega) = \mathbf{b}_{n}^T\mathbf{X}(\omega) = \mathbf{b}_{n}^T\mathbf{\phi}_r \eta_r(\omega) = \phi_{a,n} \eta_r(\omega).
				\label{eq:baseDisp}
			\end{equation}					
			Inserting Eq.~\eqref{eq:baseDisp} into Eq.~\eqref{eq:projAbs} and multiplying the first line of the latter by $1/\phi_{a,n}$, one obtains
			\begin{equation}
				\begin{array}{c} \left(-\omega^2 \left[\begin{array}{cc} \dfrac{1}{\phi_{a,n}^2} & 0 \\ 0 & m_{a,n} \end{array}\right] + j\omega \left[\begin{array}{cc} 2\zeta_r \dfrac{\omega_r}{\phi_{a,n}^2}+ c_{a,n} & -c_{a,n} \\ -c_{a,n} & c_{a,n} \end{array}\right] + \right. \\ \left. \left[\begin{array}{cc} \dfrac{\omega_r^2}{\phi_{a,n}^2} + k_{a,n} & -k_{a,n} \\ - k_{a,n} & k_{a,n} \end{array}\right]\right)\left[\begin{array}{c} U_n(\omega) \\ X_{a,n}(\omega) \end{array}\right] = \left[\begin{array}{c} \dfrac{\mathbf{\phi}_r^T \mathbf{F}(\omega)}{\phi_{a,n}} \\ 0 \end{array}\right]
				\end{array}
			\end{equation}			
			which has the same form as the equations of motions of a single-degree-of-freedom oscillator to which an absorber is attached. Warburton~\cite{Warburton1982} arrived to the same conclusion using energy considerations. 
	\section*{References}
    \bibliographystyle{elsarticle-num}
    \bibliography{biblio.bib}

 \end{document}